\documentclass[eqsecnum,showpacs,showkeys,floats,aps,nofootinbib,preprint]{revtex4}
\usepackage{bm,amsfonts, mathtools}
\usepackage{graphicx,color, wasysym}
\usepackage{textcomp}
\usepackage{amsmath,amssymb,latexsym,epsfig}

\def\de{\mathrm{d}}
\def\e{\mathrm{e}}
\def\Re{\mathrm{Re}}
\def\nn{\nonumber}

\begin{document}

\title{On evaluation of integrals involving Bessel functions}

\author{D. Babusci}
\email{danilo.babusci@lnf.infn.it}
\affiliation{INFN - Laboratori Nazionali di Frascati, via E. Fermi, 40, IT 00044 Frascati (Roma), Italy}

\author{G. Dattoli}
\email{dattoli@frascati.enea.it}
\affiliation{ENEA - Centro Ricerche Frascati, via E. Fermi, 45, IT 00044 Frascati (Roma), Italy}

\begin{abstract}
We introduce a symbolic method for the evaluation of definite integrals containing combinations of various functions, 
including exponentials, logarithm and products of Bessel functions of different types. The method we develop is naturally 
suited for the evaluation of integrals associated with specific Feynman diagrams.
\end{abstract}

\maketitle


\section{Introduction}\label{Intro}
The elementary integral
\begin{equation}
L_\nu = \int \de x\, x^{\,\nu} = \frac{x^{\,\nu + 1}}{\nu + 1} + c
\end{equation}
diverges in the limit $\nu \to - 1$, but the singularity can be removed by choosing as arbitrary constant 
\begin{equation}
c = - \frac1{\nu + 1}\,, \nn
\end{equation}
i.e., by subtracting infinities of the same order, in a way to obtain the finite result
\begin{equation}
\label{eq:log}
\lim _{\nu \to - 1}\,L_\nu = \ln x\,.
\end{equation}
This definition is very useful to evaluate integrals containing logarithmic functions. For example, its use makes a 
simple task to prove the following identity (largely used in the next sections):
\begin{align}
\label{eq:Emu}
E_\mu (a, b) &=   \int_0^\infty \de x\, x^{\,\mu} \,\e^{- a\,x^2}\,\ln (b\,x) \qquad\qquad 
(\Re\,\mu > - 1, \Re\,a > 0, b > 0) \nn \\
                      &= \frac1{4\,\sqrt{a^{\,\mu + 1}}}\,\Gamma \left(\frac{\mu + 1}2\right)\,
\left[\ln \left(\frac{b^2}{a}\right) + \psi \left(\frac{\mu + 1}2\right)\right]  
\end{align}
where $\psi$ is the digamma function. For further applications 
of the definition (\ref{eq:log}), see Ref. \cite{Babusci}.

The ``trick" of subtracting infinities of the same order is at the origin also of the definition of Macdonald function 
of integer order \cite{Spanier}
\begin{equation} 
\label{eq:Mac}
K_n (x) = \frac{\pi}2\,\lim_{\nu \to n} \frac{I_{- \nu} (x) - I_{\nu} (x)}{\sin (\nu\,\pi)} 
\end{equation}
where $I_\nu (x)$ is the hyperbolic Bessel function \cite{Spanier2} defined by the series 
\begin{equation}
\label{eq:hypBes}
I_\nu (x) = \sum_{k = 0}^\infty \frac{(x/2)^{2\,k + \nu}}{k!\,\Gamma (k + \nu + 1)}\,.
\end{equation}
In Ref. \cite{Babusci2}, the cylindrical Bessel function has been written as a product of elementary functions as 
follows
\begin{equation}
\label{eq:Bes}
J_\nu (x) = \left(\frac{x}2\,\hat{c}\right)^\nu \,\exp\left\{- \hat{c}\,\left(\frac{x}2\right)^2\right\}\,\varphi (0)
\end{equation}
where the umbral operator $\hat{c}$ is defined by its action on $\varphi (0)$ as follows:
\begin{equation}
\label{eq:cop}
\hat{c}^{\,\mu}\,\varphi (0) = \varphi (\mu) = \frac1{\Gamma (\mu + 1)}\,.
\end{equation}
As a consequence of the relationship 
\begin{equation}
\label{eq:JIrel}
J_\nu (x) = i^{- \nu}\,I_\nu (i\,x)\,,
\end{equation} 
from eq. (\ref{eq:Bes}) one has 
\begin{equation}
I_\nu (x) = \left(\frac{x}2\,\hat{c}\right)^\nu \,\exp\left\{\hat{c}\,\left(\frac{x}2\right)^2\right\}\,\varphi (0)\,,
\end{equation}
and, thus, for the Macdonald functions (\ref{eq:Mac}) the expression
\begin{align}
\label{eq:Kn}
K_n (x) &= \frac{\pi}2\,\lim_{\nu \to n} \frac{\left(\displaystyle \frac{x}2\,\hat{c}\right)^{- \nu} - 
\left(\displaystyle \frac{x}2\,\hat{c}\right)^\nu}{\sin (\nu\,\pi)}\,\exp\left\{\hat{c}\,\left(\frac{x}2\right)^2\right\}\,\varphi (0) \nn \\
              &= - \pi\, \lim_{\nu \to n} \frac{\sinh \left[\nu\,\ln \left(\displaystyle \frac{x}2\,\hat{c}\right)\right]}{\sin (\nu\,\pi)}\,
\exp\left\{\hat{c}\,\left(\frac{x}2\right)^2\right\}\,\varphi (0)\,.
\end{align}

As an example, let us consider the case of the function $K_0 (x)$. By applying l'H\^{o}pital rule, one has
\begin{align}
\label{eq:K0}
K_0 (x) &= - \ln \left(\frac{x}2\,\hat{c}\right)\, \exp\left\{\hat{c}\,\left(\frac{x}2\right)^2\right\}\,\varphi (0) \nn \\
             &= - \ln \left(\frac{x}2\right)\,I_0 (x) - \exp\left\{\hat{c}\,\left(\frac{x}2\right)^2\right\}\,(\ln \hat{c})\,\varphi (0)\,,
\end{align}
i.e.,  
\begin{equation}
K_0 (x) =  - \ln \left(\frac{x}2\right)\,I_0 (x) -  \sum_{k = 0}^\infty \frac{(x/2)^{2\,k}}{k!}\,
(\hat{c}^{\,k}\,\ln \hat{c})\,\varphi (0)\,. \nn
\end{equation}
Taking into account the definitions  (\ref{eq:log}) and (\ref{eq:cop}), we can write
\begin{align}
\label{eq:logc}
(\hat{c}^{\,\mu}\,\ln \hat{c})\,\varphi (0) = \lim_{\nu \to 0} \frac{\hat{c}^{\,\mu + \nu} - \hat{c}^{\,\mu}}{\nu}\, \varphi (0) &=
\lim_{\nu \to 0} \frac{\varphi (\mu + \nu) - \varphi (\mu)}{\nu} \nn \\
&= \varphi^\prime (\mu) = - \frac{\psi (\mu + 1)}{\Gamma (\mu + 1)}\,,
\end{align}
and, therefore
\begin{equation}
K_0 (x) = - \ln \left(\frac{x}2\right)\,I_0 (x) +  \sum_{k = 0}^\infty \frac{(x/2)^{2\,k}}{k!}\,\frac{\psi (k + 1)}{\Gamma (k + 1)}\,.  \nn 
\end{equation}
Since 
\begin{equation}
\psi (k + 1) = - \gamma + h_k \qquad\qquad  h_k = \sum_{m = 1}^k \frac1m\,, \nn
\end{equation}
(with $h_0 = 0$) where $\gamma$ is Euler's constant, we, finally, obtain
\begin{equation}
K_0 (x) = -\left[\gamma + \ln \left(\frac{x}2\right)\right]\,I_0 (x) + \sum_{k = 1}^\infty \frac{h_k}{(k!)^2}\,\left(\frac{x}2\right)^{2\,k}\,,
\end{equation}
in agreement with the expression reported in Ref. \cite{Spanier3}.

In the following sections the calculation technique outlined in this introduction will be used for the evaluation of some definite 
integrals often encountered in applications, in particular in the study of Feynman amplitudes.

\section{Integrals involving Bessel functions}\label{IntBes}
In this section, in order to fix the main results which will be exploited in the second part of the paper,  we will derive old and new 
integrals involving Bessel-type functions. The starting point will be the general procedure put forward in Ref. \cite{Babusci2} to treat the 
integrals of Bessel functions.  The technique is fairly efficient and straightforward and benefits from the fact that, as already remarked, 
Bessel functions can formally be treated as elementary functions. 

The integral ($\lambda \geq 0$) 
\begin{equation}
\label{eq:Imula}
I_{\mu, \lambda} (p) = \int_0^\infty \de x\, x^\mu\,J_\lambda (p\,x)\qquad \qquad
\left(- \Re\,\lambda - 1 < \Re\,\mu < \frac12, p > 0 \right)
\end{equation}
occurs very often in applications, for example in the evaluation of mass corrections to the large square momentum behavior of some QCD 
processes \cite{Groote}.  Taking into account the definition (\ref{eq:Bes}), we find, in a fairly straightforward way, the well known 
result (see formula \textbf{6.561.14} in Ref. \cite{GradRyz})
\begin{equation}
\label{eq:ResI}
I_{\mu, \lambda} (p) = \frac12\,\left(\frac2{p}\right)^{\mu + 1}\,
\frac{\displaystyle \Gamma\left(\frac{\lambda + \mu + 1}2\right)}{\displaystyle \Gamma\left(\frac{\lambda - \mu + 1}2\right)}\,.
\end{equation}
In the same way, it's easy to show that 
\begin{align}
\label{eq:Alam}
A_\lambda (p, b) &= \int_0^\infty \de x\, J_\lambda (p\,x)\,\ln (b\,x) \qquad\qquad  (\Re\,\lambda > - 1, p > 0, b > 0) \nn \\ 
&= \int_0^\infty \de x\, \ln (b\,x)\, \left(\frac{p\,x}2\,\hat{c}\right)^\lambda\,
\exp\left\{- \hat{c}\,\left(\frac{p\,x}2\right)^2\right\}\,\varphi (0) \nn \\
&= \left(\frac{p}2\,\hat{c}\right)^\lambda \,
E_\lambda \left(\frac{p^2}4\,\hat{c}, b\right)\,\varphi (0) 
\end{align}
and therefore, taking into account eq. (\ref{eq:Emu}), one obtains
\begin{equation}
A_\lambda (p, b) = \frac1{p}\,\left[\ln \left(\frac{2\,b}{p}\right) + \psi \left(\frac{\lambda + 1}2\right)\right]\,.
\end{equation}
Moreover, by using the expression (\ref{eq:K0}), we can write ($\Re\,\mu > - 1$)
\begin{align}
\Theta_{\mu,0} = \int_0^\infty \de x\, x^{\,\mu}\, K_0 (x) &= - \int_0^\infty \de x\, x^{\,\mu}\,\ln \left(\frac{x}2\,\hat{c}\right)\,
\exp\left\{\hat{c}\,\left(\frac{x}2\right)^2\right\}\,\varphi (0) \nn \\
&= - E_\mu \left(- \frac{\hat{c}}4, \frac{\hat{c}}2\right)\,\varphi (0) 
\end{align}
and, by using eqs. (\ref{eq:Emu}) and (\ref{eq:logc}),
\begin{equation}
\Theta_{\mu,0} = 2^{\,\mu - 1}\,\Gamma^2\left(\frac{\mu + 1}2\right)\,,
\end{equation}
that in the case $\mu = 0$ gives
\begin{equation}
\int_0^\infty \de x\, K_0 (x) = \frac{\pi}2\,. 
\end{equation}

The Neumann function (i.e., the function associated with the second solution of the cylindrical Bessel's equation)  
of integer order is defined by a procedure analogous to the one leading to Macdonald function, namely  
\begin{align}
\label{eq:Neu}
Y_n (x) &= \lim_{\nu \to n} \frac{\cos (\nu\,\pi)\,J_\nu (x) - J_{- \nu} (x)}{\sin (\nu\,\pi)} \nn \\
&= \lim_{\nu \to n} \frac{\left(\displaystyle \frac{x}2\,\hat{c}\right)^\nu\,\cos \nu\pi - 
\left(\displaystyle \frac{x}2\,\hat{c}\right)^{- \nu}}{\sin (\nu\,\pi)}\,\exp\left\{- \hat{c}\,\left(\frac{x}2\right)^2\right\}\,\varphi (0)
\end{align}
In the case $n = 0$ one gets
\begin{equation}
Y_0 (x) = \frac2{\pi}\,\ln\left(\frac{x}2\,\hat{c}\right)\,\exp\left\{- \hat{c}\,\left(\frac{x}2\right)^2\right\}\,\varphi (0)\,.
\end{equation}
By using this expression, and eq. (\ref{eq:Emu}), we can easily prove the identity 
\begin{equation}
\Upsilon_{\mu, 0} =  \int_0^\infty \de x\, x^{\,\mu}\, Y_0 (x) = \frac{2^{\,\mu}}{\pi}\,\Gamma^2\left(\frac{\mu + 1}2\right)\,
\sin\left(\mu\,\frac{\pi}2\right) \qquad\qquad \left(- 1 < \Re\,\mu < \frac12\right)
\end{equation}
that, specialized to the case $\mu = 0$, reproduces the well known result
\begin{equation}
\int_0^\infty \de x\, Y_0 (x) = 0\,.
\end{equation}

\section{Integrals involving products of Bessel functions}\label{IntProdBes}
We have already stressed that integrals of the type (\ref{eq:Imula}) appears in the solutions to many physical problems. 
In this section we will consider some generalizations of this integral associated with specific Feynman diagrams 
\cite{Groote, Mendels}. 

The product of two cylindrical Bessel functions is given by 
\begin{equation}
J_\mu (a\,x)\,J_\nu (b\,x) = \left(\frac{a\,x}2\right)^\mu\,\left(\frac{b\,x}2\right)^\nu\,\sum_{n = 0}^\infty \frac{(- 1)^n}{n!}\,
\varphi_{\mu, \nu} (n; a, b)\,\left(\frac{x}2\right)^{\,2\,n}
\end{equation}
with
\begin{equation}
\label{eq:phimn}
\varphi_{\mu, \nu} (n; a, b) = n!\,a^{\,2\,n}\,\sum_{k = 0}^n \frac{(b / a)^{\,2\,k}}
{(n - k)!\,k!\,\Gamma (n - k + \mu + 1)\,\Gamma (k + \nu + 1)}\,.
\end{equation}
Applying the symbolic method outlined before, this product can formally be rewritten as:
\begin{equation}
J_\mu (a\,x)\,J_\nu (b\,x) = \left(\frac{a\,x}2\right)^\mu\,\left(\frac{b\,x}2\right)^\nu\,
\exp\left\{- \hat{d}\,\left(\frac{x}2\right)^2\right\}\,\varphi_{\mu, \nu} (0; a, b)\,,
\end{equation}
where the operator $\hat{d}$ is defined in such a way that
\begin{equation}
\hat{d}^{\,\lambda} \,\varphi_{\mu, \nu} (0; a, b) = \varphi_{\mu, \nu} (\lambda; a, b)
\end{equation}
with $\varphi_{\mu, \nu} (\lambda; a, b)$ generalization to the noninteger case of the function introduced in 
eq. (\ref{eq:phimn}), i.e.
\begin{equation}
\varphi_{\mu, \nu} (\lambda; a, b) = \Gamma (\lambda + 1)\,a^{\,2\,\lambda}\,\sum_{k = 0}^\infty \frac{(b / a)^{\,2\,k}}
{\Gamma (\lambda - k + 1)\,k!\,\Gamma (\lambda - k + \mu + 1)\,\Gamma (k + \nu + 1)}\,.
\end{equation}
By using these results we easily get the following identity (with $a$ and $b$ unequal positive number): 
\begin{align}
\label{eq:ommn}
\Omega_{\alpha, \mu, \nu} (a, b) &= \int_0^\infty \de x\, x^\alpha\,J_\mu (a\,x)\,J_\nu (b\,x) \qquad\qquad
(\Re \left(\alpha + \mu + \nu)  >  - 1,  \Re\,\alpha < 1\right) \nn \\
&= 2^{\,\alpha}\,a^{\,\mu} \,b^{\,\nu}\,\Gamma \left(\frac{\alpha + \mu + \nu + 1}2\right)\,
\varphi_{\mu, \nu} \left(- \frac{\alpha + \mu + \nu + 1}2; a, b\right)\,, 
\end{align}
i.e., 
\begin{align}
\Omega_{\alpha, \mu, \nu} (a, b) &= \frac{\pi}2\,\left(\frac{a}2\right)^{\,\mu}\,\left(\frac{b}2\right)^{\,\nu}\,
\sec\left[(\alpha + \mu + \nu)\,\frac{\pi}2\right] \nn \\ 
& \qquad \times \,\sum_{k = 0}^\infty \frac{(b / a)^{\,2\,k}}{\Gamma \left(\displaystyle \frac{1 - \alpha - \mu - \nu}2 - k\right)\,
k!\,\Gamma \left(\displaystyle \frac{1 - \alpha + \mu - \nu}2 - k\right)\,\Gamma (k + \nu + 1)}
\end{align}
This result allows us to evaluate in a very simple way integrals involving products of different Bessel-type 
functions, as, for example,  
\begin{equation}
\Xi_{\alpha, \mu, \nu} (a, b) = \int_0^\infty \de x\, x^\alpha\,J_\mu (a\,x)\,Y_\nu (b\,x)\,,
\end{equation}
where the definition of Neumann function of noninteger order is obtained removing the limit in eq. (\ref{eq:Neu}). Taking 
into account eq. (\ref{eq:ommn}), one has
\begin{equation}
\Xi_{\alpha, \mu, \nu} (a, b) = \frac1{\sin (\nu\,\pi)}\,\left\{\cos (\nu\,\pi)\,\Omega_{\alpha, \mu, \nu} (a, b) - 
\Omega_{\alpha, \mu, -\nu} (a, b)\right\}\,.
\end{equation}

It's evident that the method here outlined can be applied, with a negligible increase in calculation complexity,  also to integrals 
involving an arbitrary number of Bessel-type functions.

\section{Concluding remarks}\label{ConRem}
In closing the paper, we apply the operatorial method described in the previous sections to the Tricomi-Bessel function 
\begin{equation}
C_\nu (x) = \frac{J_\nu (2\,\sqrt{x})}{x^{\,\nu/2}} = \sum_{k = 0}^\infty \frac{(- 1)^k\,x^k}{k!\,\Gamma (\nu + k + 1)}\,,
\end{equation}
that can be formally written as 
\begin{equation}
C_\nu (x) = \hat{c}^{\,\nu}\,\e^{- \hat{c}\,x}\,\varphi (0)
\end{equation}  
with the operator $\hat{c}$ satisfying the identity (\ref{eq:cop}). As a consequence of this equation it's easily checked that
\begin{equation}
\int_0^\infty \de x\, C_\nu (x) = \frac1{\Gamma (\nu)}\,, \qquad\qquad
\int_0^\infty \de x\, C_\nu (x^2) = \frac{\sqrt{\pi}}2\,\frac1{\Gamma \left(\nu + \frac12\right)}\,.
\end{equation}

Moreover, from Eqs. (\ref{eq:Imula}), (\ref{eq:ResI}) one has
\begin{equation}
\int_0^\infty \de x\,x^{\,\mu}\,C_\nu (x^2) = I_{\mu - \nu, \nu} (2) 
= \frac12\,\frac{\Gamma \left(\displaystyle \frac{\mu + 1}2\right)}{\Gamma \left(\nu - \displaystyle \frac{\mu - 1}2\right)}
\qquad \qquad \left(- 1 < \Re\,\mu < \frac12\right).
\end{equation}

This function, even though not explicitly mentioned as Tricomi-Bessel function, occurs in the evaluation of some Feynman 
diagrams through integrals of the type 
\begin{align}
\Psi_{\lambda, \mu, \nu} (a, b) &= \int_0^\infty \de x\,x^{\,\lambda}\,C_\mu \left(\frac{a^2\,x^2}4\right)\,K_\nu (b\,x) \nn \\
& \qquad\qquad (\Re\,(b \pm i\,a) > 0, \Re\,(2\,\mu - \lambda + 1)  >  |\Re\,\nu|)
\end{align}
where the definition of Macdonald function of noninteger is obtained not taking the limit in eq. (\ref{eq:Mac}). By using 
eqs. (\ref{eq:JIrel}) and (\ref{eq:ommn}), we get
\begin{equation}
\Psi_{\lambda, \mu, \nu} (a, b) = \frac{2^{\,\mu - 1}\,\pi}{i^{\,\nu}\,a^{\,\mu}\,\sin (\nu\,\pi)}\,
\left[(- 1)^{\,\nu}\,\Omega_{\lambda - \mu, \mu, -\nu} (a, i\,b) - \Omega_{\lambda - \mu, \mu, \nu} (a, i\,b) \right]\,.
\end{equation}
The variety of topics treated in this paper clearly shows that the method we have proposed is particularly efficient to study the problem 
of the evaluation of integrals of Bessel-type functions which very often occurs in physical application. As will be discussed in a forthcoming 
investigation, the intrinsic operational nature of the method make it ideal for an automatic implementation by a symbolic manipulator.

\end{document}